\pdfoutput=1
\documentclass[12pt]{article}
\usepackage[utf8]{inputenc}
\usepackage{color}
\usepackage[dvipsnames]{xcolor}
\usepackage{amsmath}
\usepackage{parskip}
\usepackage[showframe=false]{geometry}
\usepackage{changepage}
\usepackage{graphicx}
\usepackage[refpage]{nomencl}
\graphicspath{ {images/} }
\usepackage{amsthm}
\usepackage[compact]{titlesec}         
\titlespacing{\section}{10pt}{10pt}{10pt} 
\AtBeginDocument{
  \setlength\abovedisplayskip{10pt}
  \setlength\belowdisplayskip{10pt}}
\usepackage{tikz}
\usetikzlibrary{positioning}
\usepackage{amssymb}
\usepackage{tikz}
\usepackage{mathtools}
\usepackage{float}

\newcommand{\tpmod}[1]{{\@displayfalse\pmod{#1}}}

\usepackage[shortlabels]{enumitem}

\usepackage[english]{babel}
\usepackage{soul}
\theoremstyle{definition}
\newtheorem{theorem}{Theorem}[section]

\newtheorem{definition}[theorem]{Definition}

\newtheorem{personalEnvironment}[theorem]{$\boldsymbol{(}\mkern-3mu$}
\newcommand{\Alt}[1]{\mathrm{Alt}({#1})}
\newcommand{\Sym}[1]{\mathrm{Sym}({#1})}
\newcommand{\Dih}[1]{\mathrm{Dih}({#1})}
\newcommand{\PGL}{\mathrm{PGL}}
\newcommand{\GL}{\mathrm{GL}}
\newcommand{\Syl}{\mathrm{Syl}}
\newcommand{\PSL}{\mathrm{PSL}}
\newcommand{\SL}{\mathrm{SL}}
\renewcommand{\O}{\mathrm{O}}

\newcommand{\SO}{\mathrm{SO}}
\newcommand{\GF}{\mathrm{GF}}
\newcommand{\M}{\mathrm{M}}
\newcommand{\G}{\mathrm{G}}
\newcommand{\diag}{\mathrm{diag}}
\usepackage{graphicx}
\title{Unravelled Abstract Regular Polytopes}
\date{2020}
\author{Robert Nicolaides and Peter Rowley}
\begin{document}
\maketitle
\begin{abstract}
    This paper introduces the notion of an unravelled abstract regular polytope, and proves that $\SL_3(q) \rtimes <t>$, where $t$ is the transpose inverse automorphism of $\SL_3(q)$, possesses such polytopes for various congruences of $q$. A large number of small examples of such polytopes are given, along with extensive details of their various properties.
\end{abstract}
\section{Introduction}
Abstract polytopes trace their roots back to classical geometric objects such as the Platonic solids and more generally convex polytopes and non-convex \textquotedblleft star" polytopes in Euclidean space. During the twentieth century, with the work of Coxeter, Grünbaum, Dang and Schulte, the present day foundations of the subject of abstract polytopes evolved. Abstract regular polytopes are those abstract polytopes whose automorphism group acts regularly on its set of (maximal) flags.

Abstract regular polytopes with finitely many flags may be viewed entirely within a group theoretic framework, the basic counterpart being that of a string C-group, $G$. We will review these ideas in detail in Section 2, but for now remark that the key feature of $G$ is a generating set of involutions satisfying certain properties - we call such a set of involutions a C-string for $G$.

There is an extensive, recent, literature which explores abstract regular polytopes from the perspective of group theory. Typically, a particular group or family of groups is investigated usually with the aim of discovering if it is  the automorphism group of some abstract regular polytope. For a small selection of these see \cite{p1,p3,p4,p6,p9,p10,p11,p13,p16,p18} as well as their references. For many groups there can be a significant number of abstract regular polytopes of which the given group is the automorphism group (or, equivalently, many C-strings for that group). Nevertheless there have been successful attempts to catalogue particular abstract regular polytopes - see \cite{atlas1,atlas2}. 

The purpose of this paper is to introduce a certain type of abstract regular polytope which we now define. 

\begin{definition}
Suppose $G$ is a finite group with $\{t_1,\ldots,t_n\}$ a C-string for $G$ of rank $n$.
\begin{enumerate}[$(i)$]
    \item For $N \trianglelefteq G,$ set $\overline{G} = G/N$ and use $\overline{g}$ for the image of $g$ in $\overline{G}$, $g \in G$. If either $|\{\overline{t_1}, \ldots, \overline{t_n}\}| < n$ or, $|\{\overline{t_1}, \ldots, \overline{t_n}\}| = n$ and $\{\overline{t_1}, \ldots, \overline{t_n}\}$ is not a C-string for $\overline{G},$ we say that  $\{t_1,\ldots,t_n\}$ is an $N$-unravelled C-string for $G$.
    \item If $\{t_1,\ldots,t_n\}$ is an $N$-unravelled 
    C-string for all non-trivial proper normal 
    subgroups $N$ of $G$, we call $\{t_1,\ldots,t_n\}$ 
    an unravelled C-string for $G$, and refer to the 
    associated abstract regular polytope as being 
    unravelled. 
\end{enumerate}
\end{definition}

Clearly, if $G$ is a simple group, then all its C-strings will be unravelled (and all its abstract regular polytopes likewise). The motivation for looking at unravelled C-strings was to focus upon a smaller set of C-strings, that were, in some sense, intrinsic to $G$. Hence with this in mind, we need to concentrate upon non-simple groups. Even with only a few non-trivial proper normal subgroups, the unravelled C-strings often stand out from the crowd. For example $G = \SL_3(7) \rtimes \langle t \rangle \sim 3^{\boldsymbol{\cdot}} L_3(7):2,$ where $t$ acts upon $\SL_3(7)$ as the transpose inverse automorphism, has 3256 abstract regular polytopes only one of which (of rank 4) is unravelled. We note here that we shall be using the \textsc{Atlas} \cite{TheATLAS} conventions to describe the shape of groups.

In geometric terms we are considering the rank $n$ abstract regular polytopes that do not cover any others of the same rank. One could naturally extend this definition in a variety of ways, for example, dropping the restriction in rank or regularity. There exists much literature examining instances of this including 
\cite{AllPolytopess},
\cite{Moreonquotientpolytopes},
\cite{Simplertestsforsemisparsesubgroups},
\cite{Covers},
\cite{Quotientsofpolytopes}. We restrict ourselves to only the unravelled case.

Another example, related to a sporadic group is 
$3^{\boldsymbol{\cdot}}\M_{22}:2$ which has 727 abstract regular polytopes 
among which five are unravelled (and all of rank 4).

The unravelled example for $3^{\boldsymbol{\cdot}} L_3(7):2$ already mentioned appears as a specific case in our first theorem.

\begin{theorem}\label{thm 1.2} Suppose that $q$ is a prime power and $G = \SL_3(q) \rtimes <t>$ where $t$ acts upon $\SL_3(q)$ as the transpose inverse automorphism. Assume that 
\begin{enumerate}[$(i)$]
    \item $6 | q-1$;
    \item there exists $\lambda, \mu \in \GF(q) $ such that $2\lambda^2 = 1$ and $2\mu^2 - \lambda^2 = 0$; and 
    \item at least one of $-3^{-1}+(3^{-2}-1)^{{1}/{2}}$ and $-3^{-1}-(3^{-2}-1)^{{1}/{2}}$ has order $q+1$ in $\GF(q^2)^*$. 
\end{enumerate}
Then $G$ possesses an unravelled rank 4 C-string with Schläfli symbol $[4,q+1,4].$
\end{theorem}
There are infinitely many $q$ satisfying $(i)$ and $(ii)$ of Theorem \ref{thm 1.2} (for example, taking $q=p$, a prime with  $p \equiv 1\pmod 3$ and $p\equiv 7 \pmod 8$ gives infinitely many $q$ by Dirichlet's Theorem). However we do not know if there are infinitely many $q$ satisfying all three conditions in the theorem. Of the 157 primes $p$ less than or equal to 10000 with $p \equiv 1 \pmod 3$ and $p\equiv 7 \pmod 8$, 20 of them do not satisfy $(iii)$ (and they are $
199, 343, 919, 1039, 1063, 2239, 3079, 3919, 4423, 4759, 4783, 5167, 6967, 7039, 7759,\\ 7879, 8287, 8887, 9511, 9679)$.

A \textsc{Magma}\cite{Magma} calculation shows there are no unravelled C-strings for $\SL_3(13):2 \sim 3^{\boldsymbol{\cdot}}L_3(13):2$ with Schläfli symbols $[4,14,4]$. So the requirement of $p \equiv 7 \pmod 8$ in Theorem \ref{thm 1.2} is necessary. However, we have the following companion result to Theorem \ref{thm 1.2}.

\begin{theorem}\label{thm 1.3} Let $p$ be a prime with $p \equiv 1 \pmod 3$ and $p \equiv 5 \pmod 8.$ Then $G=\SL_3(p) \rtimes \langle t \rangle,$ where $t$ is the transpose inverse automorphism of $\SL_3(p)$, has an unravelled rank 4 C-string with Schläfli symbol $[4,p,4]$. 

\end{theorem}

In Nicolaides and Rowley \cite{BnUNravlledPaper} two further families of unravelled C-strings are uncovered, both associated with Coxeter groups of type $B_n$. For one family they all have rank 4, but the other has unbounded rank.

The proofs of Theorems \ref{thm 1.2} and \ref{thm 1.3} occupy Sections 3 and 4, respectively. Moreover, these proofs are constructive. The involutions of the C-strings being given as explicit $6 \times 6$ matrices. As already mentioned Section 2 reviews notation and concepts relevant to this paper. Our final section is a compendium of calculations, done with the aid of \textsc{Magma}\cite{Magma}. Table \ref{table1} is a census of the number, up to isomorphism, of C-strings for a variety of groups. This table also records how many of the C-strings are unravelled. Among these groups we have the Coxeter groups of types $B_3,B_4,B_5,B_6,B_7,B_8,D_3,D_4,D_5,D_6,D_7,D_8$, joined by a number of groups which are unusual in some respect. We single out for mention the non-split extensions $3^{\boldsymbol{\cdot}}\Sym{6}$ and $3^{\boldsymbol{\cdot}}\Sym{7}$. Each of these groups has a unique unravelled C-string of rank 4. The remainder of this section highlights various properties of these along with other unravelled C-strings.

\section{Preliminaries}
Although we shall work exclusively in the group theory context, we say a few words about abstract regular polytopes.

An abstract regular polytope $\mathcal{P}$ is a ranked partially ordered poset, whose elements are usually called faces, and satisfies four axioms. First it's also assumed that $\mathcal{P}$ has a smallest and largest face. We assume the ranks are precisely $I = \{1,\ldots,n\}$, so $\mathcal{P}$ has rank $n$, with the smallest and largest faces conventionally given ranks $0$ and $n+1$ respectively. A flag (or chamber) of $\mathcal{P}$ is a maximal totally ordered chain in $\mathcal{P}$. Two flags $F_1$ and $F_2$ of $\mathcal{P}$ are said to be $i$-adjacent where $i \in I$ if they differ in exactly one face of rank (or type) $i$. The chamber graph of $\mathcal{P}$ has the flags of $\mathcal{P}$ as its vertices with two flags adjacent in the chamber graph if they are $i$-adjacent for some $i \in I$. The other three axioms that $\mathcal{P}$ must fulfill are that all flags of $\mathcal{P}$ should have the same number of faces, it be strongly connected and satisfy the diamond property. For $\mathcal{P}$ to be strongly connected each of its sections (included $\mathcal{P}$) must be connected. The definition of sections, and further background, may be found in Sections 2B and 2E of McMullen, Schulte \cite{ARP}. The diamond property requires that for every flag $F$ and $i \in I$, there is precisely one other flag in $\mathcal{P}$ which is $i$-adjacent to $F$.

An abstract regular polytope $\mathcal{P}$ is called regular if its automorphism group, $Aut(\mathcal{P})$ acts transitively on the set of flags of $\mathcal{P}$. This leads, via $Aut(\mathcal{P})$, to being able to translate abstract regular polytopes entirely into group theoretic data, see McMullen, Schulte \cite{ARP} again, Section 2E for details. That group entity is a string C-group. Suppose $G$ is a group, and let $\{t_1,\ldots,t_n\}$ be a set of involutions in $G$. Set $I = \{1,\ldots,n\}$, and for $J\subseteq I$ define $G_J = \langle t_j | j \in I\rangle$, setting $G_J =1$ when $J = \emptyset$. Sometimes we write $G_J$ as $G_{j_1\dotsm j_k}$ when $J=\{j_1,\ldots,j_k\}$. We say $\{t_1,\ldots,t_n\}$ is a C-string for $G$ if 
\begin{enumerate}[$(i)$]
    \item $G = \langle t_1,\ldots,t_n \rangle$;
    \item $t_i$ and $t_j$ commute whenever $|i-j| \ge 2$; 
    and
    \item for $J,K \subseteq I$, $G_J \cap G_K = G_{J \cap K}$.
\end{enumerate}
We note that $(iii)$ is usually referred to as the intersection property. A group possessing a C-string (and a given group may have many) is called a string C-group. The Schläfli symbol of a C-string $\{t_1,\ldots,t_n\}$ is the sequence $[\tau_{12},\tau_{23},\ldots,\tau_{n-1n}]$ where $\tau_{jj+1}$ is the order of $t_jt_{j+1}.$
In Sections 3 and 4 we will be investigating C-strings in the group $G = \SL_3(q) \ltimes
 \langle t \rangle $ where $q$ is some prime power with $3 \mid q-1$ and $t$ is the transpose inverse automorphism of $\SL_3(q)$. We close this section establishing some relevant notation. Put $H = \SL_3(q)$ and let $U$ be the natural 3-dimensional $\GF(q)H$-module. Set $V = U \oplus U^*$, where $U^*$ is the dual of $U$. Choosing a basis for $U$ and a dual basis for $U^*$ (viewing $U$ and $U^*$ as subspaces of $V$) we may take $t$ to be   $ t = \left(\begin{array}{c|c}
& I_3 \\
\hline
I_3 & \\
\end{array}\right). $
We note that $G$ has two conjugacy classes of involutions, namely $t^G$ and $s^G$ where $s \in G' = H$. These classes may be easily distinguished as $\dim C_V(t) = 3$ whereas $\dim C_V(s) = 2$. Also, since $3 \mid q-1$, $G$ has shape $3^{\boldsymbol{\cdot}}L_3(q):2.$
Our group theoretic notation is standard as given, for example, in \cite{suzuki1982group}. Additionally, we use $\Sym{n}$, and $\Alt{n}$ to denote, respectively, the symmetric and alternating groups of degree $n$.

\section{C-strings with Schläfli symbol $[4,q+1,4]$}
In this section we prove Theorem \ref{thm 1.2} in a series of steps. We use the set up given at the end of Section 2.
Since $6 \mid q-1$, we may select $\rho \in \GF(q)^\ast $ such that $\rho$ has multiplicative order 6. Further, we have $\lambda, \mu \in \GF(q)$ for which $2\lambda^2=1$ and $2\mu^2-\lambda^2=0$. We now introduce five other elements of $\GF(q)$.

\par\quad\par \begin{personalEnvironment}
$\mkern-7mu \boldsymbol{)}$  \label{noteworthyElements}
\begin{align*}
    \alpha &= (\mu^{-2}-1)^{-1}\\
    \beta &= 2(\mu^{-2}-1)^{-1}\lambda\mu^{-1}\\
    \xi &= \rho^2+(1-\rho^2)2^{-1}\\
    \eta &= (1-\rho^2)2^{-1}\\
    \tau &= \rho^4
\end{align*}
\end{personalEnvironment}

Observe that $\alpha = 3^{-1}$. From $2\mu^2 = \lambda^2 = 2^{-1}$ we get $2^{-1}\mu^{-2} = 2$, and so $\mu^{-2} = 4.$ Therefore $\alpha = (\mu^{-2} -1)^{-1} = 3^{-1}$. Also, since $\beta = 2\alpha\lambda\mu^{-1}$, $$\beta^2 = 4\alpha^2\lambda^2\mu^{-2} = 8\alpha^2.$$

Hence $\alpha^2 + \beta^2 = 9\alpha^2 = 9(3^{-1})^2 = 1.$ Thus $\alpha^2 + \beta^2 = 1.$

Using these elements we now define our C-string, $\{t_1,t_2,t_3,t_4\}$. We shall show that $\{t_1,t_2,t_3,t_4\}$ is an unravelled C-string for $G$ where the $t_i$ are specified as follows.

\begin{personalEnvironment}$\boldsymbol{\mkern-7mu)}$\label{generators}
\begin{align*}
   t_1 &= \left(\begin{array}{c c c|c c c }
& & & \mu & \phantom{-}\lambda &\phantom{-}\mu\\
&\phantom{-}\makebox(0,0){\text{\huge0}} & & \lambda & \phantom{-}0 &-\lambda\\
& & & \mu & -\lambda &\phantom{-}\mu\\
\hline
 \mu & \phantom{-}\lambda &\phantom{-}\mu& & &\\
\lambda & \phantom{-}0 &-\lambda& &\phantom{-}\makebox(0,0){\text{\huge0}} &\\
\mu & -\lambda &\phantom{-}\mu& & &\\
\end{array}\right) \\
   t_2 &= \left(\begin{array}{c c c|c c c }
-1& & & & & \\
& 1& & &\makebox(0,0){\text{\huge0}} & \\
& &  -1& & &\\
\hline
& & &-1& &  \\
& \makebox(0,0){\text{\huge0}}& && 1 & \\
& & && & -1\\
\end{array}\right) = \diag(-1,1,-1,-1,1,-1) \\
   t_3 &= \left(\begin{array}{c c c|c c c }
\phantom{-}\alpha& \phantom{-}\beta& \phantom{-}0& & & \\
\phantom{-}\beta& -\alpha&\phantom{-}0 & &\makebox(0,0){\text{\huge0}} & \\
\phantom{-}0&\phantom{-}0 &  -1& & &\\
\hline
& & &\phantom{-}\alpha&\phantom{-}\beta &  \phantom{-}0\\
& \phantom{-}\makebox(0,0){\text{\huge0}}& &\phantom{-}\beta& -\alpha &\phantom{-}0 \\
& & &\phantom{-}0&\phantom{-}0 & -1\\
\end{array}\right) \\
   t_4 &= \left(\begin{array}{c c c|c c c }
& & & \xi & 0 &\eta\\
& \makebox(0,0){\text{\huge0}}& & 0 & \tau &0\\
& & & \eta &0 &\xi\\
\hline
 \xi\rho^{-2} & 0 &\eta\rho& & &\\
0 & \tau\rho^{-2} &0& &\makebox(0,0){\text{\huge0}} &\\
\eta\rho &0 &\xi\rho^{-2}& & &\\
\end{array}\right) \\
\end{align*}
\end{personalEnvironment}

\begin{personalEnvironment}$\boldsymbol{\mkern-7mu )}$\label{involutionsLemma}
For $i=1,2,3,4$, $t_i$ are involutions with $t_1,t_4 \in t^G$ and $t_2,t_3 \in s^G$.
\end{personalEnvironment}
\par The diagonal blocks of $t_2$ and $t_3$ are easily seen to be involutions, and so $t_2$ and $t_3$ are involutions. Since
\begin{equation*}
  \left(\begin{array}{c c c }
\mu & \lambda &\mu\\
 \lambda & 0 &-\lambda\\
\mu & -\lambda &\mu\\
\end{array}\right)^2 \\ = 
\left(\begin{array}{c c c }
2\mu^2+\lambda^2 & 0 &2\mu^2-\lambda^2\\
 0 & 2\lambda^2 &0\\
2\mu^2-\lambda^2 &0 &2\mu^2+\lambda^2\\
\end{array}\right) \\
\end{equation*}
 the conditions on $\mu$ and $\lambda$ imply that $t_1$ is  an involution.
 \par Moving onto $t_4$, we look at the product
 \begin{equation*}
  \left(\begin{array}{c c c }
 \xi & 0 &\eta\\
 0 & \tau &0\\
 \eta &0 &\xi\\
\end{array}\right)
  \left(\begin{array}{c c c }
 \xi\rho^{-2} & 0 &\eta\rho\\
0 & \tau\rho^{-2} &0\\
\eta\rho &0 &\xi\rho^{-2}\\
\end{array}\right)
\\ = 
  \left(\begin{array}{c c c }
 \xi^2\rho^{-2}+\eta^2\rho & 0 &\xi\eta\rho+\eta\xi\rho^{-2}\\
0 & \tau^2\rho^{-2} &0\\
\eta\xi\rho^{-2}+\xi\eta\rho & 0 &\eta^2\rho+\xi^2\rho^{-2}\\
\end{array}\right) = A.
\end{equation*}
Note that $\rho^3$ has multiplicative order 2, and so $\rho^3=-1$.
Now 
\begin{align*}
    \eta\xi\rho^{-2}+\xi\eta\rho&=\eta\xi\rho^{-2}(1+\rho^3)\\
    &= \eta\xi\rho^{-2}(1+-1)=0,
\end{align*}
and using (\ref{noteworthyElements}) we have 
\begin{align*}
    \tau^2\rho^{-2}=\rho^8\rho^{-2}=\rho^6=1
\end{align*}
Again, from (\ref{noteworthyElements})
\begin{align*}
   \xi&=\rho^2+\eta && \\
   \xi^2&=\rho^4+2\rho^2\eta+\eta^2 && \\
   \xi^2\rho^{-2}&=\rho^2+2\eta+\eta^2\rho^{-2} && \\
   \xi^2\rho^{-2}+\eta^2\rho&=\rho^2+2\eta+\eta^2\rho^{-2}+ \eta^2\rho&& \\
   &=\rho^2+(1-\rho^2)+\eta^2\rho^{-2} +\eta^2\rho&& \\
\end{align*} 
as $2\eta=1-\rho^2$. Then, as $\eta^2\rho^{-2}+\eta^2\rho = \eta^2\rho^{-2}(1+\rho^3)=0$, we get $$\xi^2\rho^{-2}+\eta^2\rho = 1. $$
Hence $A=I_3$, whence $t_4$ is also an involution. Since $\dim C_V(t_i)=3$ for $i=1,4$ and $\dim C_V(t_i)=2$ for $i=2,3$,  (\ref{involutionsLemma}) is proved.

\par\quad\par\begin{personalEnvironment}$\boldsymbol{\mkern-7mu )}$\label{lemma4}
\begin{align*}
    C_G(t) &= \langle t \rangle \times C_H(t) \cong 2 \times \SO_3(q) \cong 2 \times \PGL_2(q)
\end{align*}
Because $t$ acts by inverse conjugation on $H$, $C_H(t)$ consists of all orthogonal matrices of determinant 1. The well-known isomorphism $\SO_3(q) \cong \PGL_2(q)$ (see \cite{taylor1992geometry}) now gives   (\ref{lemma4}).
\end{personalEnvironment}
\par We define 
$$   r = \left(\begin{array}{c c c|c c c }
& & & \rho & 0 &0\\
& \makebox(0,0){\text{\huge0}}& & 0 & \rho &0\\
& & & 0 &0 &\rho^{-2}\\
\hline
 \rho^{-1} & 0 &0& & &\\
0 & \rho^{-1} &0& &\makebox(0,0){\text{\huge0}} &\\
0 &0 &\rho^{2}& & &\\
\end{array}\right) \\. $$
Observe that $r \in t^G$ and so $C_G(r) \cong 2 \times \PGL_2(q).$

\par\quad\par\begin{personalEnvironment}$\boldsymbol{\mkern-7mu )}$\label{lemma5}
$tr= \diag(\rho^{-1},\rho^{-1},\rho^{2},\rho^{},\rho^{},\rho^{-2}) \in H$ has order 6 and $(tr)^2 \in Z(H)$. Further, $C_G(t) \cap C_G(r) \le C_G(tr) = C_H(tr)  \cong \GL_2(q)$. 

Since $[G:H]=2,$ we have $tr \in H$ and, as $\rho$ has multiplicative order 6, $tr$ has order 6 with $(tr)^2 \in Z(H).$ Thus $C_G(tr)=C_H(tr)=C_H((tr)^3) \cong \GL_2(q).$
\end{personalEnvironment}

\par\quad\par\begin{personalEnvironment}$\boldsymbol{\mkern-7mu )}$\label{lemma6}
We have $t_1,t_2,t_3 \in C_G(t)$ and $t_2,t_3,t_4 \in C_G(r)$.
\end{personalEnvironment}

It  is straightforward to check (\ref{lemma6}), though for $t_4r=rt_4$ we use the fact that $\rho^2=\rho^{-4}$.

\par\quad\par\begin{personalEnvironment}$\boldsymbol{\mkern-7mu )}$\label{lemma7}
$C_G(t) \cap C_G(r) \cong \Dih{2(q+\epsilon)}$ where $\epsilon = \pm 1.$
\end{personalEnvironment}
First we observe that $C_G(t) \cap C_G(r) = C_{C_G(tr)}(t).$ Since $C_G(tr)=C_H(tr)\cong \GL_2(q)$ by   (\ref{lemma5}) and $t$ acts by transpose inverse upon $C_H(tr)$, $C_{C_G(tr)}(t) \cong \O_2^\epsilon(q)$ (the 2-dimensional orthogonal group of type $\epsilon$). Since $\O_2^\epsilon(q) \cong \Dih{2(q-\epsilon)}$, (see \cite{taylor1992geometry}), we have (\ref{lemma7}).

\par\quad\par\begin{personalEnvironment}$\boldsymbol{\mkern-7mu )}$\label{lemma8}
The order of $t_1t_2$ is 4.
\end{personalEnvironment}
We have $t_1t_2 = \left(\begin{array}{c | c }
&A\\
\hline
A&
\end{array}\right)$ where $A= \left(\begin{array}{c c c }
-\mu&\phantom{-}\lambda&-\mu\\
-\lambda&\phantom{-}0&\phantom{-}\lambda\\
-\mu&-\lambda&-\mu\\
\end{array}\right)$. Now $A^2 =  \left(\begin{array}{c c c }
2\mu^2-\lambda^2&0&2\mu^2+\lambda^2\\
0&-2\lambda^2&0\\
2\mu^2+\lambda^2&0&2\mu^2-\lambda^2\\
\end{array}\right)$ and hence $A^2 =  \left(\begin{array}{c c c }
0&\phantom{-}0&\phantom{-}1\\
0&-1&\phantom{-}0\\
1&\phantom{-}0&\phantom{-}0\\
\end{array}\right)$.
Therefore $t_1t_4$ has order 4.

\par\quad\par\begin{personalEnvironment}$\boldsymbol{\mkern-7mu )}$\label{lemma9}
$t_1t_3 = t_3t_1$.
\end{personalEnvironment}
Let $A= \left(\begin{array}{c c c }
\mu&\phantom{-}\lambda&\phantom{-}\mu\\
\lambda&\phantom{-}0&-\lambda\\
\mu&-\lambda&\phantom{-}\mu\\
\end{array}\right)$
and $B= \left(\begin{array}{c c c }
\alpha&\phantom{-}\beta&\phantom{-}0\\
\beta&-\alpha&\phantom{-}0\\
0&\phantom{-}0&-1\\
\end{array}\right).$
Then $t_1t_3 = t_3t_1$ provided $AB=BA$. Now 
\begin{align*}
    AB &= \left(\begin{array}{c c c }
\mu\alpha + \lambda\beta&\mu\beta - \alpha\lambda&-\mu\\
\lambda\alpha&\lambda\beta&\phantom{-}\lambda\\
\mu\alpha - \lambda\beta&\mu\beta + \alpha\lambda&-\mu\\
\end{array}\right) &&\text{ and }\\
BA&= \left(\begin{array}{c c c }
\alpha\mu + \beta\lambda&\alpha\lambda&\alpha\mu - \beta\lambda\\
\beta\mu - \alpha\lambda&\beta\lambda&\beta\mu + \alpha\lambda\\
-\mu&\phantom{\alpha}\lambda&-\mu\\
\end{array}\right).
\end{align*}
So we need to know that 
\begin{align*}
    \alpha\lambda &= \mu\beta-\alpha\lambda,\\
    -\mu &= \alpha\mu-\beta\lambda \quad \text{and}\\
    \lambda &= \beta\mu+\alpha\lambda.\\
\end{align*}
Since $\mu\beta = \mu2(\mu^{-2}-1)^{-1}\lambda\mu^{-1}=2(\mu^{-2}-1)^{-1}\lambda = 2\alpha\lambda,$ we have $\alpha\lambda = \mu\beta-\alpha\lambda$. From $\lambda\beta = \lambda2(\mu^{-2}-1)^{-1}\lambda\mu^{-1}=(\mu^{-2}-1)^{-1}\mu^{-1}=\alpha\mu^{-1}$,
we get 
\begin{align*}
    \mu\alpha-\lambda\beta &= \mu\alpha - \alpha\mu^{-1}\\
                           &= \mu\alpha(1-\mu^{-2})\\
                           &= \mu(\mu^{-2}-1)^{-1}(1-\mu^{-2})\\
                           &= -\mu.
\end{align*}
Finally we show $\lambda=\beta\mu+\alpha\lambda$. Using $\beta=2\alpha\mu^{-1}$, we have 
\begin{align*}
    \mu\beta + \alpha\lambda &= 2\alpha\lambda+\alpha\lambda\\
    &= 3\alpha\lambda \\
    &= 3(\mu^{-2}-1)^{-1}\lambda\\
    &=3.3^{-1}\lambda = \lambda,
\end{align*}
as $4\mu^2=1$ implies $\mu^{-2}-1=3.$ Hence (\ref{lemma9}) holds.

\par\quad\par\begin{personalEnvironment}$\boldsymbol{\mkern-7mu )}$\label{lemma10}
The order of $t_2t_3$ is $q+1$ and $C_G(t) \cap C_G(r)  = \langle t_2,t_3\rangle$
\end{personalEnvironment}

We use that 

\begin{align*}
    t_2t_3 &= \left(\begin{array}{c | c }
X&\\
\hline
&X
\end{array}\right)    \text{where  } X = \left(\begin{array}{c c c }
-\alpha & -\beta&0\\
\beta & -\alpha&0\\
0&0&1\\
\end{array}\right).
\end{align*}

Hence the order of $t_2t_3$ is the same as the order of $Y$ where $Y = \left(\begin{array}{c c  }
-\alpha & -\beta\\
\beta & -\alpha\\
\end{array}\right).$ Recalling that $\alpha^2 + \beta^2 = 1$, the characteristic polynomial of $Y$ is $$ x^2 +2\alpha x +1.$$
Therefore the eigenvalues of $Y$ are $-\alpha \pm (\alpha^2-1)^{1/2} = -3^{-1}\pm(3^{-2}-1)^{1/2}.$ If these two eigenvalues are equal, then $2(\alpha^2 - 1)^{1/2} = 0$ which implies the impossible $\alpha^2 = 1$. So the two eigenvalues of $Y$ are different. Consequently $Y$ is diagonalizable in $\GL_2(q^2)$ and hence, by assumption $(iii)$ of Theorem \ref{thm 1.2}, $Y$ has order $q+1$. Hence, using (\ref{lemma6}) and (\ref{lemma7}), we obtain $C_G(t) \cap C_G(r) = \langle t_2,t_3 \rangle.$

\par\quad\par\begin{personalEnvironment}$\boldsymbol{\mkern-7mu )}$\label{lemma11.5}
$[t_2,t_4] = 1$
\end{personalEnvironment}
Since $t_2$ is a diagonal matrix with 1 and -1 as its only diagonal entries, a matrix commutes with $t_2$ if and only if its of the form
$$\left(\begin{array}{c c c |c c c }
*&0&*&*&0&*\\
0&*&0&0&*&0\\
*&0&*&*&0&*\\
\hline
*&0&*&*&0&*\\
0&*&0&0&*&0\\
*&0&*&*&0&*\\
\end{array}\right),$$
and $t_4$ is of this form.

\par\quad\par\begin{personalEnvironment}$\boldsymbol{\mkern-7mu )}$\label{lemma12}
$[t_1,t_4] = 1$
\end{personalEnvironment}
Writing $t_1 = \left(\begin{array}{c | c }
&A\\
\hline
A&
\end{array}\right)$ and $t_4 = \left(\begin{array}{c | c }
&C\\
\hline
D&
\end{array}\right)$, (\ref{lemma12}) will hold if we show that $AD=CA$ and $AC=DA$. Calculating gives
\begin{align*}
    AD &= \left(\begin{array}{c c c }
\mu\xi\rho^{-2}+\mu\eta\rho&\lambda\tau\rho^{-2}&\mu\eta\rho+\mu\xi\rho^{-2}\\
\lambda\xi\rho^{-2}-\lambda\eta\rho^{}&0&\lambda\eta\rho-\lambda\xi\rho^{-2}\\
\mu\xi\rho^{-2}+\mu\eta\rho &-\lambda\tau\rho^{-2}&\mu\eta\rho+\mu\xi\rho^{-2}\\
\end{array}\right) \quad\text{and}\\
 CA &= \left(\begin{array}{c c c }
\xi\mu + \eta\mu&\xi\lambda-\eta\lambda&\xi\mu + \mu\eta\\
\tau\lambda&0&-\tau\lambda\\
\eta\mu+\xi\mu&\eta\lambda-\xi\lambda&\eta\mu+\xi\mu\\
\end{array}\right).
\end{align*}
Therefore $AD=CA$ holds provided
\begin{align*}
    \mu\xi\rho^{-2}+ \mu\eta\rho &= \xi\mu+\eta\mu,\\
    \lambda\xi\rho^{-2}- \lambda\eta\rho &= \tau\lambda \quad \text{and}\\
    \lambda\tau\rho^{-2} &= \xi\lambda-\eta\lambda.\\
\end{align*}
Since $\lambda \ne 0$ and $\mu \ne 0$ this is equivalent to showing that
\begin{align*}
    \xi\rho^{-2}+\eta\rho &= \xi+\eta, \\
     \xi\rho^{-2}-\eta\rho &= \tau \quad \text{and} \\
     \tau\rho^{-2} &= \xi-\eta.
\end{align*}
First we observe that $\xi = \rho^2 + \eta$, and recall that $\rho^3=-1$. Hence 
\begin{align*}
    \xi+\eta &= \rho^2+2\eta \\
    &= \rho^2+2(1-\rho^2)2^{-1}\\
    &=\rho^2+1-\rho^2 = 1.
\end{align*}
While 
\begin{align*}
    \xi\rho^{-2}+\eta\rho &= (\rho^2+\eta)\rho^{-2} + \eta\rho\\
    &=1+\eta\rho^{-2} + \eta\rho\\
    &=1+\eta\rho^{-2}(1+\rho^3)\\
    &=1+\eta\rho^{-2}(1-1)=1.
\end{align*}
Next,
\begin{align*}
    \xi\rho^{-2}-\eta\rho &= (\rho^2+\eta)\rho^{-2}- \eta\rho\\
    &= 1 + \eta\rho^{-2}-\eta\rho\\
    &=\rho^4(\rho^2+\eta-\eta\rho^{-3})\\
    &=\rho^4(\rho^2+2\eta),
\end{align*}
and substituting for $\eta$ yields
\begin{align*}
    \xi\rho^{-2} - \eta\rho&= \rho^4(\rho^2+2(1-\rho^2)2^{-1})\\
    &=\rho^4=\tau.
\end{align*}
Since $\xi -\eta= \rho^2+\eta-\eta=\rho^2=\rho^4\rho^{-2}=\tau\rho^{-2},$we have shown that $AD=CA.$ Similar considerations verify that $AC=DA$, whence (\ref{lemma12}) holds.

\par\quad\par\begin{personalEnvironment}$\boldsymbol{\mkern-7mu )}$\label{lemma12.5}
$t_3t_4$ has order 4. 
\end{personalEnvironment}

Let 
\begin{align*}
    X &= \left(\begin{array}{c c c }
\alpha & \beta & 0\\
\beta & -\alpha & 0\\
0 & 0 & -1
\end{array}\right),\\
    A &= \left(\begin{array}{c c c }
\xi& 0 & \eta\\
0 & \tau & 0\\
\eta & 0 & \xi 
\end{array}\right) and\\
    B &= \left(\begin{array}{c c c }
\xi\rho^{-2}& 0 & \eta\rho\\
0 & \tau\rho^{-2} & 0\\
\eta\rho & 0 & \xi\rho^{-2} 
\end{array}\right).
\end{align*}

To show that $t_3t_4$ has order 4 we verify that $(t_3t_4)^2$ is an involution. Now

\begin{align*}
    (t_3t_4)^2 &= \left(\begin{array}{ c | c }
XAXB& \\
\hline
&XBXA\\
\end{array}\right).
\end{align*}
We will see in a moment that the $(3,2)^{th}$-entry of XBXA is non-zero, so $(t_3t_4)^2 \ne 1$. Thus recalling that $X=X^{-1}$ and $A^{-1}= B$, we must show 
\begin{align*}
    XAXB &= (XAXB)^{-1} = AXBX \text{ and }\\
    XBXA &= (XBXA)^{-1} = BXAX.
\end{align*}
Observe that $XBXA=BXAX$ implies
    $$A(XBXA)B = A(BXAX)B,$$ giving $AXBX = XAXB.$ Hence it suffices to show that $XBXA = BXAX.$
    
    We calculate that

\begin{align*}
    BXAX &= \left(\begin{array}{c c c }
\alpha^2\xi^2\rho^{-2}-\alpha\eta^2\rho+\beta^2\xi\tau\rho^{
-2} & \alpha\beta\xi^2\rho^{-2}-\beta\eta^2\rho-\alpha\beta
\xi\tau\rho^{-2} & \eta\xi\rho-\alpha\xi\eta\rho^{-2}\\
\alpha\beta\xi\tau\rho^{-2}-\alpha\beta\tau^2\rho^{-2} & 
\beta^2\xi\tau\rho^{-2}+\alpha^2\tau^2\rho^{-2} & 
-\beta\eta\tau\rho^{-2}\\
\alpha^2\xi\eta\rho-\alpha\xi\eta\rho^{-2}+\beta^2\tau\eta
\rho & \alpha\beta\xi\eta\rho-\beta\xi\eta\rho^{-2}-\alpha
\beta\tau\eta\rho & \eta^2\rho^{-2}-\alpha\xi^2\rho\\
\end{array}\right) \text{and}   \\
    XBXA &= \left(\begin{array}{c c c }
\alpha\xi^2\alpha\rho^{-2}-\alpha\eta^2\rho+\beta^2\tau{\rho}^{-2}\xi & \alpha\xi{\rho}^{-2}\beta\tau - \beta\tau^2\rho^{-2}\alpha & \alpha^2\xi\rho^{-2}\eta - \alpha\eta\rho\xi+\beta^2\tau\rho^{-2}\eta\\
\beta\xi^2\alpha\rho^{-2}-\beta\eta^2\rho-\alpha\tau\rho^{-2}\beta\xi & \beta^2\xi\rho^{-2}\tau + \alpha^2\tau^2\rho^{-2} & \beta\xi\rho^{-2}\alpha\eta-\beta\eta\rho\xi-\alpha\tau\rho^{-2}\beta\eta\\
-\eta\rho\alpha\xi + \xi\rho^{-2}\eta & \eta\rho\beta\tau & - \eta^2 \rho\alpha+\xi^2\rho^{-2}
\end{array}\right).   
\end{align*}

First we note that $XBXA$ and $BXAX$ have the same diagonal entries. For the $(2,1)^{th}$-co-ordinate of $XBXA$ and $BXAX$ we require $$\alpha\beta\xi\tau\rho^{-2}-\alpha\beta\tau^2\rho^{-2} = \beta\xi^2\alpha\rho^{-2}-\beta\eta^2\rho-\alpha\tau\rho^{-2}\beta\xi.  $$
Multiplying through by $\beta\rho^2$ this is equivalent to
$$\tau\xi\alpha - \alpha\tau^2 = \xi^2\alpha+\eta^2 - \alpha\tau\xi, $$ using $\beta^3 = -1$. Since $-2\xi = \tau$, this is equivalent to $$-2\alpha\tau^2 = \xi^2\alpha+\eta^2. $$
Substituting for $\xi,\eta$ and $\alpha = 3^{-1}$ reduces this to
$$0 = 1 + \rho^2 +\rho^4, $$ which holds. Therefore $XBXA$ and $BXAX$ have the same $(2,1)^{th}$-co-ordinate. Similarly we may check all the off-diagonal entries of $XBXA$ and $BXAX$ are equal. Therefore $XBXA = BXAX$ and hence (\ref{lemma12.5}) holds.

\par\quad\par\begin{personalEnvironment}$\boldsymbol{\mkern-7mu )}$\label{lemma11}
$\langle t_1,t_2,t_3 \rangle = C_G(t)$ and $\langle t_2,t_3,t_4 \rangle = C_G(r)$.
\end{personalEnvironment}
Since $t_1,t_2 \in C_G(t)$ with $t_2 \in C_H(t) \trianglelefteq C_G(t)$, $[t_1,t_2] \in C_H(t) \cong \PGL_2(q)$. Now, from   (\ref{lemma8}), $$[t_1,t_2] = (t_1t_2)^2 =  \left(\begin{array}{c c c |c c c }
&&1&&&\\
&-1&&&&\\
1&&&&&\\
\hline
&&&&&1\\
&&&&-1&\\
&&&1&&\\
\end{array}\right).$$
A quick calculation reveals that $[t_1,t_2] \notin C_G(r)$, and so $[t_1,t_2] \notin C_G(t)\cap C_G(r).$ By (\ref{lemma7}) and \cite{suzuki1982group}, $C_G(t) \cap C_G(r)$ is a maximal subgroup of $C_H(t)$, whence, as $t_1 \notin C_H(t),$ we infer that $\langle t_1,t_2,t_3 \rangle = C_G(t)$. Similar considerations show that $\langle t_2,t_3,t_4\rangle = C_G(r).$
\par\quad\par\begin{personalEnvironment}$\boldsymbol{\mkern-7mu )}$\label{lemma13}
$\{t_1,t_2,t_3,t_4 \}$ is a C-string for $G$ with Schläfli symbol $[4,q+1,4].$
\end{personalEnvironment}
This comes from combining the fact that $C_H(t)$ is a maximal subgroup of $H$ (see \cite{Mitchell}) with (\ref{lemma8}), (\ref{lemma9}), (\ref{lemma10}), (\ref{lemma11.5}), (\ref{lemma12}) and (\ref{lemma12.5}).

\par\quad\par\begin{personalEnvironment}$\boldsymbol{\mkern-7mu )}$\label{lemma14}
$\{t_1,t_2,t_3,t_4 \}$ is an unravelled C-string of $G$. \par 
\end{personalEnvironment}
The only  non-trivial proper normal subgroups of $G$ are $H$,$Z(H)$. Since $[G:H]=2$, we only need show $\{t_1,t_2,t_3,t_4\}$ is $Z(H)$-unravelled. 
Put $G_{123}=\langle t_1,t_2,t_3\rangle$, $G_{234}=\langle t_2,t_3,t_4\rangle$ and $\overline{G}=G/Z(H)$. 
Since $Z(H)= \langle \diag(\rho^{-2},\rho^{-2},\rho^{-2},\rho^{2},\rho^{2},\rho^{2}) \rangle$, we see that $|\{\overline{t_1},\overline{t_2},\overline{t_3},\overline{t_4}\}|=4$. 
Also $\langle \overline{t_2},\overline{t_3}\rangle \cong C_G(t) \cap C_G(r) \cong \Dih{2(q+1)}$ as $\langle {t_2},{t_3}\rangle \cap Z(H) = 1$. 
Now $\overline{G}_{123}=C_{\overline{G}}(\overline{t})$ and $\overline{G}_{234}=C_{\overline{G}}(\overline{r}),$ as the orders of $t$ and $r$ are coprime to $|Z(H)|$. 
From   (\ref{lemma5}) $\overline{t}\overline{r}=\overline{r}\overline{t}$ has order 2. 
That is $\overline{t}$ and $\overline{r}$ commute. So $\overline{t},\overline{r} \in \overline{G}_{123} \cap \overline{G}_{234}$ and therefore $\overline{G}_{123} \cap \overline{G}_{234} \gneqq
 \langle \overline{t_2},\overline{t_3} \rangle$. 
 Consequently, the intersection property fails for $\{\overline{t_1},\overline{t_2},\overline{t_3},\overline{t_4}\}$.
Together (\ref{lemma13}) and (\ref{lemma14}) prove Theorem \ref{thm 1.2}.

\section{C-strings with Schläfli symbol $[4,p,4]$}
As mentioned in Section 2, this section is concerned with proving Theorem \ref{thm 1.3}. Again we employ the notational set up described at the end of Section 2. Now here $G = \SL_3(p) \rtimes \langle t \rangle$ where $p$ is a prime such that $p \equiv 1 \pmod 3$ and $p \equiv 5 \pmod 8$. Because $p \equiv 1 \pmod 3$ we may choose, and keep fixed, $\rho \in \GF(p)$ of multiplicative order 3. Further, $p \equiv 5 \pmod 8$ means we may choose $\iota \in \GF(p)$, also now to be fixed, such that $\iota^2 = -1$. Set $\alpha = \sqrt{(1+\rho^2)^{-1}}$, again making a choice from the (at most) two possibilities. Now we define a slew of elements in $\GF(p)$.
\begin{personalEnvironment}\label{(4.1)}$\boldsymbol{\mkern-7mu )}$
\end{personalEnvironment}
\par 
\begin{align*}
     \lambda &= \alpha(\iota+1)(-1+\rho-\iota\rho^2)\\
    \epsilon &= -\iota\lambda\\
    \beta &= -2^{-1}\lambda^2\iota\\
    \gamma &= 2^{-1}\lambda^2-1\\
    \delta &= -1 -2^{-1}\lambda^2\\
    \mu &= 1 - \rho.\\
\end{align*}
Note that $\lambda \ne 0$ and $\lambda^2 = -\epsilon^2.$ Also recall that $1 + \rho+\rho^2 = 0$ and so $\alpha^2 = -\rho^2.$ Hence $\alpha \ne 0$.

The elements in (\ref{(4.1)}) appear as entries in $\{t_1,t_2,t_3,t_4\}$, elements of $G$, which we now define.

\begin{personalEnvironment}\label{(4.2)}$\boldsymbol{\mkern-7mu )}$

\begin{align*}
   t_1 &= \left(\begin{array}{c c c|c c c }
& & & \phantom{-}0 & \phantom{-}\alpha &-\alpha\rho\\
&\phantom{-}\makebox(0,0){\text{\huge0}} & & \phantom{-}\alpha & \phantom{-}\rho &\phantom{-}1\\
& & & -\alpha\rho & \phantom{-}1 &\phantom{-}\rho^2\\
\hline
\phantom{-}0 & \phantom{-}\alpha &-\alpha\rho& & & \\
\phantom{-}\alpha & \phantom{-}\rho &\phantom{-}1&&\phantom{-}\makebox(0,0){\text{\huge0}} & \\
 -\alpha\rho & \phantom{-}1 &\phantom{-}\rho^2& & &\\
\end{array}\right) \\
   t_2 &= \left(\begin{array}{c c c|c c c }
1& & & & & \\
& -1& & &\makebox(0,0){\text{\huge0}} & \\
& &  -1& & &\\
\hline
& & &1& &  \\
& \makebox(0,0){\text{\huge0}}& && -1 & \\
& & && & -1\\
\end{array}\right) = \diag(1,-1,-1,1,-1,-1) \\
   t_3 &= \left(\begin{array}{c c c|c c c }
1& \lambda& \epsilon& & & \\
\lambda&\gamma&\beta & &\makebox(0,0){\text{\huge0}} & \\
\epsilon&\beta &  \delta& & &\\
\hline
& & &1& \lambda& \epsilon \\
& \makebox(0,0){\text{\huge0}} &&\lambda&\gamma&\beta  \\
& & &\epsilon&\beta &  \delta\\
\end{array}\right) \\
   t_4 &= \left(\begin{array}{c c c|c c c }
& & & -\rho & \phantom{-}0 &\phantom{-}0\\
& \makebox(0,0){\text{\huge0}}& & \phantom{-}0  &\phantom{-}0 & \phantom{-}\rho\\
& & & \phantom{-}0 &\phantom{-}\rho &-\mu\rho^2\\
\hline
 -\rho^{2} & 0 &0& & &\\
\phantom{-}0 & \mu &\rho^2& &\makebox(0,0){\text{\huge0}} &\\
\phantom{-}0 &\rho^{2}&0& & &\\
\end{array}\right).
\end{align*}
\end{personalEnvironment}
In order to define a further element in $t^G$, we introduce more elements in $\GF(p)$.
\begin{personalEnvironment}\label{(4.3)}$\boldsymbol{\mkern-7mu )}$
\begin{align*}
    a &= 2(2\rho^2+(1-\rho)\iota)^{-1}\\
    x &= -2^{-1}a\rho(1-\rho)\\
    y &= -x\rho\\
    b &= a^{-1}(\rho^2 + x^2)\\
    c &= a^{-1}(1+x^2\rho)\\
    d &= a\rho\\
\end{align*}
Observe that $a\ne0$, so $b$ and $c$ are well-defined. Now set 
\begin{align*}
    r &= \left(\begin{array}{c c c|c c c }
& & & \rho & 0 &0\\
& \makebox(0,0){\text{\huge0}}& & 0 & a &x\\
& & & 0 &x &b\\
\hline
 \rho^{2} & 0 &0& & &\\
0 & c &y& &\makebox(0,0){\text{\huge0}} &\\
0 &y &d& & &\\
\end{array}\right). \\
\end{align*}

\end{personalEnvironment}
\begin{personalEnvironment}\label{(4.4)}$\boldsymbol{\mkern-7mu )}$
\begin{enumerate}[$(i)$]
    \item $t_1,t_2,t_3,t_4$ and $r$ are involutions.
    \item $t_1,t_4,r \in t^G$ and $t_2,t_3 \in s^G.$
\end{enumerate}
To show that $t_1$ is an involution, we must verify that $X^2 = I_3$ where $ X = \left(\begin{array}{c c c}
0&\alpha&-\alpha\rho\\
\alpha&\rho&1\\
\alpha\rho&1&\rho^2\\
\end{array}\right).$

Now $X^2 = \left(\begin{array}{c c c}
\alpha^2+\alpha^2\rho & 0 &0\\
0 & \alpha^2 + \rho^2 +1 & -\alpha^2\rho+\rho+\rho^2\\
0 & -\alpha^2\rho+\rho+\rho^2& \alpha^2\rho+1+\rho^4
\end{array}\right)$
and using $\alpha^2 = -\rho^2,$ we see $X^2 = I_3.$
Similarly, using (\ref{(4.1)}), we may show $t_3$ is an involution. While it is straightforward to check that $t_2$ and $t_4$ are involutions, for $r$ it suffices, using (\ref{(4.3)}), to show that 
$$\left(\begin{array}{c  c}
a & x\\
x & b\\
\end{array}\right)^{-1}= \left(\begin{array}{c  c}
c & y\\
y & d\\
\end{array}\right), $$
so proving $(i)$. Since, by calculation, $\dim C_V(t_1) = \dim C_V(t_4) = \dim C_V(r) = 3$ and $\dim C_V(t_2) = \dim C_V(t_3) = 2,$ we have part $(ii)$.
\end{personalEnvironment}

\begin{personalEnvironment}\label{(4.5)}$\boldsymbol{\mkern-7mu )}$
$t_1t_3 = t_3t_1,$ $t_1t_4 = t_4t_1$ and $t_2t_4 = t_4t_2$.

Checking $t_1t_4 = t_4t_1$ uses $\mu = 1 - \rho$ whereas $t_1t_3 = t_3t_1$ requires the definitions of $\lambda, \epsilon, \beta,\gamma$ and $\delta$. That $t_2t_4 = t_4t_2$ is easily seen.
\end{personalEnvironment}

\begin{personalEnvironment}\label{(4.6)}$\boldsymbol{\mkern-7mu )}$
\begin{enumerate}[$(i)$]
    \item $t_1t_2$ and $t_3t_4$ both have order 4.
    \item $t_2t_3$ has order $p$.
\end{enumerate}

Part $(i)$ can be checked following the same strategy as in (\ref{lemma12.5}).

Now $t_2t_3 = \left(\begin{array}{c | c}
X &\\
\hline
& X\\
\end{array}\right)$ where $X = \left(\begin{array}{c c c}
1 & \lambda & \epsilon\\
-\lambda & -\gamma & -\beta\\
-\epsilon & -\beta & -\delta\\
\end{array}\right).$

We demonstrate that $X$ has order $p$, from which $(ii)$ will follow. Consider $X$ acting on the 3-dimensional vector space $U$, setting $U_1 = C_U(X)$ and letting $U_2$ be the inverse image of $C_{U/U_1}(X)$ in $U$. For $(u,v,w) \in U,  (u,v,w) \in  U_1$ if and only if 
\begin{align*}
   u - \lambda v - \epsilon w &= u\\
    \lambda u - \gamma v - \beta w &= v\\
    \epsilon u - \beta v - \delta w &= w\\
\end{align*} 
The first equation gives $v = - \lambda^{-1}\epsilon w = - \lambda ^{-1}. - \iota \lambda w = \iota w,$ and then the second yields $$ \lambda u = (\gamma \iota + \beta + \iota)w = 0, $$
using the definitions of $\gamma$ and $\beta$. Since $\lambda \ne 0, u = 0.$ Thus $U_1 = \{(0,\iota w, w ) | w \in \GF(p)\}$. Similar calculations show that $U_2 = \{ (u,\iota w,w) | u,w \in \GF(p)\}$. Now $(0,0,1)X - (0,0,1) = (-\epsilon,-\beta,-\delta-1) \in U_2$, as $\iota(-\delta-1) =-\beta.$ Hence as $(0,0,1) \notin U_2,$ $X$ acts nilpotently on $U$, whence $X$ has $p$-power order. Since Sylow $p$-subgroups of $\SL_3(p)$ have exponent $p$ and $X \ne I_3$, $X$ has order $p$. This completes the proof of (\ref{(4.6)}).
\end{personalEnvironment}

Let $g_0 = \left(\begin{array}{c c c|c c c }
& & & 1 & 0 &\phantom{-}0\\
& \makebox(0,0){\text{\huge0}}& & 0 & 0 &-1\\
& & & 0 &1 &\phantom{-}0\\
\hline
 1 & 0 &\phantom{-}0& & &\\
0 & 0 &-1& &\makebox(0,0){\text{\huge0}} &\\
0 &1 &\phantom{-}0& & &\\
\end{array}\right)$ and
$z = \diag(\rho,\rho,\rho,\rho^2,\rho^2,\rho^2).$ Note that $z \in Z(H)$, and straightforward calculation gives
\begin{personalEnvironment}\label{(4.7)}$\boldsymbol{\mkern-7mu )}$
\begin{enumerate}[(i)]
    \item $g_0 \in C_G(t)$ and $zg_0 \in C_G(r).$
    \item $g_0^2 = t_2 = (zg_0)^2$.
\end{enumerate}
\end{personalEnvironment}
Set $L_{123} = G_{123} \cap C_H(t)'$ and $L_{234} = G_{234} \cap C_H(t)'.$ Note that $L_{123} \cong \PSL_2(p) \cong L_{234}.$

\begin{personalEnvironment}\label{(4.8)}$\boldsymbol{\mkern-7mu )}$
\begin{enumerate}[(i)]
    \item $C_G(t) \ge G_{123}$ and $C_G(r) \ge G_{234}.$
    \item $G_{123} = \langle t_1 \rangle L_{123}$ and $G_{234} = \langle t_4 \rangle L_{234}.$
    \item $G_{123} \cong \PGL(2,p) \cong G_{234.}$
\end{enumerate}
\end{personalEnvironment}
First, calculation reveals that $t_1, t_2$ and $t_3$ commute with $t$ and $t_2,t_3$ and $t_4$ commute with $r$, so part $(i)$ holds.

Observe that, as $C_G(t)=\langle t\rangle \times C_H(t)$ with $C_H(t) \equiv \PGL(2,p)$, $C_G(t) = \langle t \rangle \times C_H(t)$ with $C_H(t) \cong \PGL(2,p), C_G(t)/L_{123}$ is elementary abelian of order 4, (\ref{(4.7)}) implies that $t_2 \in L_{123}.$ Clearly we also have $t_2t_3 \in L_{123},$ so $G_{23} = \langle t_2, t_3 \rangle \le L_{123}.$ Since by (\ref{(4.6)})$(ii)$, $\Dih{2p} \cong G_{23}$ is a maximal subgroup of $L_{123} \cong \PSL(2,p)$ and $t_1$ does not normalize $G_{23},$ $G_{123} = \langle t_1 \rangle L_{123}$. A similar argument establishes $G_{234} = \langle t_4 \rangle L_{234}.$

Since $p \equiv 5 \pmod{8}$, (\ref{(4.6)})(i) implies that $\langle t_1, t_2 \rangle \in \Syl_2 G_{123}.$ Hence $t \notin G_{123}$ and so, by $(ii)$, $G_{123} \cong \PGL(2,p).$ Likewise we have $G_{234} \cong \PGL(2,p)$, so proving (\ref{(4.8)}).

\begin{personalEnvironment}\label{(*4.9*)}$\boldsymbol{\mkern-7mu )}$
$G = \langle t_1,t_2,t_3,t_4\rangle$.
\end{personalEnvironment}
Put $\overline{G} = G/Z(H).$ Then $\overline{H} \equiv \PSL_3(p)$ and $\overline{G}_{123}$ contains a subgroup isomorphic to $\PSL_2(p)$ by (\ref{(4.7)}). Since $\overline{C_G(t)}$ is the only maximal subgroup of $\overline{G}$ containing $\overline{G}_{123}$ and $\overline{t}_4 \notin \overline{C_G(t)}$, $\overline{G} = \langle \overline{G}_{123}, \overline{t}_4 \rangle.$ Now $H$ being a non-split central extension this then implies (\ref{(*4.9*)}).

\begin{personalEnvironment}\label{(*4.10*)}$\boldsymbol{\mkern-7mu )}$
$G_{23} = G_{123} \cap G_{234} = C_G(t) \cap C_G(r)$.
\end{personalEnvironment}
From (\ref{(4.7)}) $G_{23} \le G_{123} \cap G_{234} \le C_G(t) \cap C_G(r)$. Now 

$$tr = \left(\begin{array}{c c c|c c c}
\rho^2&0&0&&&\\
0&c&y&&\makebox(0,0){\text{\huge0}}&\\
0&y&d&&&\\
\hline
&&&\rho&0&0\\
&\makebox(0,0){\text{\huge0}}&&0&a&x\\
&&&0&x&b\\
\end{array}\right). $$
Let $g=ztr$ (recall that $z = \diag(\rho,\rho,\rho,\rho^2,\rho^2,\rho^2)$. Then 
$$g =  \left(\begin{array}{c c c|c c c}
1&0&0&&&\\
0&\rho c&\rho y&&\makebox(0,0){\text{\huge0}}&\\
0&\rho y&\rho d&&&\\
\hline
&&&1&0&0\\
&\makebox(0,0){\text{\huge0}}&&0&\rho^2a&\rho^2x\\
&&&0&\rho^2x&\rho^2b\\
\end{array}\right).$$

Investigating the action of $g$ on $V$ we discover that $g$ acts nilpotently on $V$, and therefore $g$ has order $p$. Hence $tr=z^{-1}g$ has order $3p$ with $\langle z \rangle \le \langle tr \rangle.$ Consequently $C_G(tr) \le C_G(z) = H$. So $C_G(tr) = C_H(g)$. Since $G_{23} \le C_G(t) \cap C_G(r) \le C_G(tr), C_G(tr)$ has even order by (\ref{(4.6)})$(i)$. Thus from centralizers of $p$-elements in $\SL_3(p)$ we have $C_G(tr) = C_H(g) \sim p^3 : (p-1).$ Let $P \in \Syl_p C_H(g).$ Then $P \trianglelefteq C_H(g).$ Also $t$ acts upon $C_H(g)/P\cong p-1$. If $t$ centralizes $C_H(g)/P$, then $C_H(g) = C_{C_H(g)}(t)P.$ Now $\langle t_2t_3\rangle \le C_H(t)$ and from $C_H(t) \cong \PGL_2(p)$ we have $N_{C_H(t)}(\langle t_2t_3\rangle) \sim p:p-1,$ so $C_{C_H(g)}(t)$ normalizes $\langle t_2t_3\rangle$ which contradicts the structure of $C_H(g).$ Therefore $t$ does not centralize $C_H(g)/P.$ Since $C_H(g)/P$ is a cyclic group, $t$ must act by inverting which implies $C_{C_G(tr)}(t)$ has order dividing $2p^3$. But the largest power of $p$ dividing $|\PGL_2(p)|$ is $p$ and so $|C_{C_G(tr)}(t)|=2p.$ Now we infer that $C_G(t) \cap C_G(r) = C_{C_G(tr)}(t) = G_{23}.$

\begin{personalEnvironment}\label{(*4.11*)}$\boldsymbol{\mkern-7mu )}$
$\{t_1,t_2,t_3,t_3\}$ is an unravelled C-string for $G$ with Schläfli symbol $[4,p,4].$
\end{personalEnvironment}

Combining (\ref{(4.4)})$(i),$ (\ref{(4.6)}), (\ref{(*4.9*)}) and (\ref{(*4.10*)}) gives that $\{t_1,t_2,t_3,t_3\}$ is a C-string with Schläfli symbol $[4,p,4].$ We now show it is unravelled. 

Since $L_{123} \cong \PSL(2,p)$ and, by assumption $p \equiv 5 \pmod{8},$ the Sylow 2-subgroup of $L_{123}$ are elementary abelian. In particular, $L_{123}$ contains no elements of order 4. Hence, if $h$ is an element of $G_{123}$ of order 4, $G_{123} = \langle h \rangle L_{123}$. As a consequence any $G_{123}$-conjugate of $h$ is $L_{123}$-conjugate. By (\ref{(4.8)})$(iii)$ $G_{123} \cong \PGL(2,p)$ and so, as its Sylow 2-subgroups are isomorphic to $\Dih{8},$ has only one $G_{123}-$conjugacy class of elements of order 4. Now
\begin{align*}
    t_1t_2 &= \left(\begin{array}{c|c }
        0 & \ast \\
        \hline
        \ast & 0
    \end{array}\right).
\end{align*}
Thus we conclude, as $L_{123} \le H,$ that all order 4 elements of $G_{123}$ must have this shape. From \ref{(4.7)}(i) $g_0 \in C_G(t)$ and, since $C_G(t) = \langle t \rangle G_{123},$ either $g_0$ or $tg_0$ are in $G_{123}$. But $tg_0$ has shape 
$\left(\begin{array}{c|c }
    \ast & 0 \\
    \hline
    0 & \ast
\end{array}\right),$ whence we deduce that $g_0 \in G_{123}.$ Because of (\ref{(4.7)}), a similar argument yields that $zg_0 \in G_{234}.$ Let $\overline{G} = G/Z(H).$ Then, as $ z \in Z(H),$ we have $$\overline{g_0} = \overline{zg_0} \in \overline{G}_{123} \cap \overline{G}_{234}, $$
but $\overline{g_0} \notin \overline{G}_{23} = \langle \overline{t}_2, \overline{t}_3 \rangle$ as $\overline{g_0}$ has order $4$. Thus $\{t_1,t_2,t_3,t_4\}$ is an unravelled C-string, so proving (\ref{(*4.11*)}).

\section{Some small unravelled C-strings}

This section is a pot pourri of calculations, obtained with the aid of \textsc{Magma}\cite{Magma}, focussing mainly on C-strings of groups which display some kind of exceptional behaviour. 

Among the alternating groups, $\Alt{6}$ and $\Alt{7}$ stand out by virtue of being the only ones having Schur multiplier divisible by 3 (see \cite{suzuki1982group}). This means we may construct groups $G$ of shape $3^{\boldsymbol{\cdot}}\Sym{6}$ and $3^{\boldsymbol{\cdot}}\Sym{7}$ where $G'$ is isomorphic to, respectively, the non-split central extensions $3^{\boldsymbol{\cdot}}\Alt{6}$  and $3^{\boldsymbol{\cdot}}\Alt{7}$.

Another kind of unusual behaviour is having a smaller dimension matrix representation than expected . Such a phenomenon, effectively via the isomorphism $\Alt{8} \equiv \GL_4(2),$ gives rise to groups of shape $2^4:\Alt{5}\equiv \M_{20}$ ($\M_{20}$ is the \textquotedblleft Mathieu" group of degree 20), $2^4:\Sym{6},$ $2^4:\Alt{7}$ and ${2^4}:\Alt{8}$. And a further example is the non-split extension ${2^4} ^{\boldsymbol{\cdot}}\Alt{8}$ (see miscellaneous groups in \cite{wilson1999atlas}).

Table \ref{table1} below is a census of C-strings for the given groups - $B_n$ and $D_n$ as usual denoting the Coxeter groups of type $B$ and $D$ of rank $n$. In a column of the table $i(j)[k]$ indicates that $i$ is the total number of C-strings (up to isomorphism) possessing the property listed for the column, of which $j$ are self-dual and $k$ the number of the $i$ C-strings which are unravelled. As a general comment we observe that $3^{\boldsymbol{\cdot}}\G_2(3):2$ with 625 C-strings is bereft of any unravelled C-strings while $3^{\boldsymbol{\cdot}}\M_{22}:2$ and $2^4:\Sym{6}$ are both well endowed. While  ${2^4}^{\boldsymbol{\cdot}}\Alt{8}$ and $2^4:\Alt{5}$ have no C-strings at all!

\par

\begin{table}[H]
\begin{adjustwidth}{-1.5cm}{}
\begin{tabular}{c|c | c| c| c| c| c| c }
    
    Group & Total & rank 3 & rank 4 & rank 5 & rank 6 & rank 7 & rank 8\\
    \hline
     $3.S_6$ &11(3)[1]&3(1)[0] &8(2)[1] &0 &0  &0 &0 \\
     $3.S_7$ &167(5)[1]&142(4)[0] &23(1)[1] &2 &0 &0  &0 \\
     $3.\PSL_3(7):2$ &3256(48)[1]&3240(44)[0] & 16(4)[1] & 0& 0&0  &0 \\
     $3.\PSL_3(13):2$ &38594(174)[1]&38534(166)[0] & 60(8)[1] & 0& 0&0  &0  \\
     $3.\M_{22}:2$ &727(13)[5]&550(10)[0] &177(3)[5] &0 &0 &0  &0 \\
     $3.\G_2(3):2$ &725(25)[0]&705(25)[0] & 20(0)[0] & 0& 0&0  &0  \\
     $2^4:S_6$ &22(2)[11]&6(0)[0] & 8(0)[4] & 8(2)[7]& 0&0  &0  \\
     \hline
     $B3$ & 8(0)[0] & 8(0)[0] &0  &0 &0 &0 &0 \\
     $B4$ & 14(2)[0] & 6(2)[0] &8(0)[0]  &0 &0 &0 &0  \\
     $B5$ & 165(0)[0] & 63(0)[0] &88(0)[0]  &14(0)[0] &0 &0 &0 \\
     $B6$ & 130(0)[0] & 24(0)[0] &76(0)[0]  &20(0)[0] &10(0)[0] &0 &0  \\
     $B7$ &2965(21)[14] & 1031(21)[0] &1428(0)[10]  &400(0)[4] &84(0)[0] &22(0)[0] &0   \\
     $B8$ &3051(33)[38] &1020(32)[0] & 1494(0)[32] &304(0)[8]  &192(0)[0] &27(1)[0] &14(0)[0]   \\
     \hline 
     $D3$ & 3(1)[3] & 3(1)[3] &0  &0 &0 &0 &0 \\
     $D4$ & 0 & 0 &0  &0 &0 &0 &0  \\
     $D5$ & 39(1)[16] & 21(1)[0] &16(0)[14]  &2(0)[2] &0 &0 &0  \\
     $D6$ & 132(0)[2] & 24(0)[0] &48(0)[2]  &60(0)[0] &0 &0 &0  \\
     $D7$ &628(16)[210] & 348(16)[0] &226(0)[166]  &42(0)[36] &10(0)[6] &2(0)[2] &0  \\
     $D8$ &3537(27)[24] & 887(19)[0] &1598(8)[14]  &826(0)[10] &172(0)[0] &54(0)[0] &0  \\
     \hline
\end{tabular}  
    \caption{Number of C-strings}
    \label{table1}
\end{adjustwidth}
\end{table}

For $\gamma$ a chamber in an abstract regular polytope and $i \in \mathbb{N}$, $\Delta_i(\gamma)$ consists of chambers which are distance $i$ from $\gamma$ in the chamber graph.

\begin{personalEnvironment}\label{(5.1)}$\boldsymbol{\mkern-7mu )}$
$G\sim 3^{\boldsymbol{\cdot}}\Sym{6}$
\end{personalEnvironment}
Up to isomorphsim $G$ has precisely one unravelled C-string. It has rank 4 and its Schläfli symbol is $[4,5,4]$ with Betti numbers $[1,18,135,13,18,1]$. The disc structure from an arbitary chamber $\gamma$ of the associated abstract regular polytope is 
\begin{figure}[H]
\centering
\begin{adjustwidth}{-0.5cm}{}
\begin{tabular}{c|c c c c c c c c c c c c c c c c c}
    
    i&1&2&3&4&5&6&7&8&9&10&11&12&13&14&15&16\\
    \hline
    $|\Delta_i(\gamma)|$  & 4 & 9 & 18
& 34 & 61 & 108 & 162 & 218 & 303 & 358 & 373 & 276 & 154 & 70 & 9
& 2  \\
     
\end{tabular}
\end{adjustwidth}
\end{figure}
So the diameter of this graph is 16. If we take $\gamma$ to correspond to the identity element of $G$, then the two chambers in $\Delta_{16}(\gamma)$ correspond to the two non-trivial elements in $Z'(G).$ Moreover, the three chambers $\{\gamma\} \cup \Delta_{16}(\gamma)$ are mutually at distance 16. We briefly examine some other properties of this C-string, \tikzstyle{int}=[draw, fill=black!20, minimum size=2em]
\tikzstyle{init} = [pin edge={to-,thin,black}]
\begin{tikzpicture}[baseline={(0,-0.5)}]

\draw (0,-0.5) -- (1,-0.5);
\draw (1,-0.5) -- (2,-0.5);
\draw (2,-0.5) -- (3,-0.5);
\filldraw [black] (0,-0.5) circle (3pt);
\filldraw [black] (1,-0.5) circle (3pt);
\filldraw [black] (2,-0.5) circle (3pt);
\filldraw [black] (3,-0.5) circle (3pt);
\node (p1) at (0.5,-0.2) {$4$};
\node (p2) at (1.5,-0.2) {$5$};
\node (p3) at (2.5,-0.2) {$4$};

\node (t1) at (0,-1) {$t_1$};
\node (t2) at (1,-1) {$t_2$};
\node (t3) at (2,-1) {$t_3$};
\node (t4) at (3,-1) {$t_4$};
\end{tikzpicture}. The subgroups $G_{123}$ and $G_{234}$ of $G$ (corresponding to the so-called vertices and facets) are both isomorphic to $\Sym{5}$. Yet $G_{123}$ and $G_{234}$ are not conjugate in $G$. According to Hartley's Atlas \cite{atlas1} there is only one abstract regular polytope of $\Sym{5}$ with Schläfli type $[4,5]$ and $[5,4]$. They are both locally spherical, non-orientable compact quotients of hyperbolic space.

\begin{personalEnvironment}\label{(5.2)}$\boldsymbol{\mkern-7mu )}$
$G\sim 3^{\boldsymbol{\cdot}}\Sym{7}$
\end{personalEnvironment}
Just as in (\ref{(5.1)}), $G$ has up to isomorphism exactly one unravelled C-string. Also it has rank 4, while its Schläfli symbol is $[4,6,4]$ and Betti numbers are $[1,63,945,945,63,1].$ For $\gamma$ a chamber of the associated abstract regular polytope, the disc sizes of the chamber graph are 
\begin{table}[H]
\centering
\begin{tabular}{c | c c c c c c c c c c c }
    i&1&2&3&4&5&6&7&8&9&10&11\\
    \hline
    $|\Delta_i(\gamma)|$   &4 & 9 & 18
& 34 & 62 & 113 & 204 & 366 & 601 & 963 &1454  \\ 
\end{tabular}
\begin{tabular}{c c c c c c c c c c c c}
12&13&14&15&16&17&18&19&20&21&22\\
\hline 
 2036 & 2562 & 2696&  2005 &1219& 514 & 188 & 57 & 10 & 4 & 1  \\
\end{tabular}
\end{table}
So we have a unique chamber at distance 22 from $\gamma$. If $\gamma$ corresponds to the identity element of $G$, then this chamber will correspond to an involution of $G$. We remark that both vertices and facets of this polytope have automorphism group isomorphic to $\mathbb{Z}_2 \times \Sym{5}$ and are non-orientable, spherical hyperbolic quotients named as $\{4,6\} \ast 240a$ polytope in Hartley's Atlas \cite{atlas1}.

\begin{personalEnvironment}\label{(5.3)}$\boldsymbol{\mkern-7mu )}$
$G\sim 3^{\boldsymbol{\cdot}}\M_{22}:2$
\end{personalEnvironment}
Here there are five unravelled C-strings all of rank 4, with details given in Table \ref{table2}. 
\begin{table}[H]
\centering
\begin{tabular}{c|c}
    
    Schläfli Symbol & Betti Numbers\\
    \hline
    
    $[4,5,4]$&$ [ 1, 2016, 166320, 166320, 8316, 1 ]$\\
    $[4,5,4]$&$ [ 1, 8316, 166320, 166320, 2016, 1 ]$\\
    $[4,6,4]$&$[ 1, 693, 166320, 166320, 693, 1 ]$\\
    $[4,6,4]$&$[ 1, 693, 166320, 166320, 6930, 1 ]$\\
    $[4,6,4]$&$[ 1, 6930, 166320, 166320, 693, 1 ]$\\
    \hline
\end{tabular} 
    \caption{Unravelled C-strings for $3^{\boldsymbol{\cdot}}\M_{22}:2$}
    \label{table2}
\end{table}
We note that the five polytopes in Table \ref{table2} consists of a dual pair of $[4,5,4]$ polytopes and a dual pair of $[4,6,4]$ polytopes and one self-dual $[4,6,4]$ polytope.

\begin{personalEnvironment}\label{(5.4)}$\boldsymbol{\mkern-7mu )}$
$G\sim 2^4:\Sym{6}$
\end{personalEnvironment}
Here in Table \ref{table3} we find eleven unravelled C-strings, four of which have rank 4 and the remainder rank 5.
\begin{table}[H]
\centering
 \begin{tabular}{ c | c }

Schläfli Symbol & Betti Numbers\\
\hline
    $[6,6,4 ]$	   &$[ 1, 60, 720, 480, 16, 1 ]$\\
    $[ 4, 6, 6 ]$	   &$[ 1, 16, 480, 720, 60, 1 ]$\\
    $[ 6, 5, 4 ]$	  &  $[ 1, 72, 720, 480, 16, 1 ]$\\
    $[ 4, 5, 6 ]$	   &$[ 1, 16, 480, 720, 72, 1 ]$\\
    $[ 4, 4, 6, 3 ]$ & $[ 1, 16, 120, 240, 90, 6, 1 ]$\\
    $[ 3, 6, 4, 4 ]$&    $[ 1, 6, 90, 240, 120, 16, 1 ]$\\
    $[ 4, 4, 4, 3 ]$ &  $[ 1, 16, 120, 240, 90, 10, 1 ]$\\
    $[ 3, 4, 4, 4 ]$   &$[ 1, 10, 90, 240, 120, 16, 1 ]$\\
     $[ 3, 6, 4, 3 ]$&	$[ 1, 6, 120, 320, 120, 16, 1 ]$\\
     
    $[ 3, 4, 6, 3 ]$   &$[ 1, 16, 120, 320, 120, 6, 1 ]$\\
    $[ 3, 4, 4, 3 ]$   &$[ 1, 16, 120, 320, 120, 16, 1 ]$\\   
    
\hline
\end{tabular}
\caption{Unravelled C-strings for $2^4:\Sym{6}$}
\label{table3}
\end{table}

Only two of the eleven, namely those with symbols $[4,5,6]$ and $[6,5,4]$, decrease in rank when quotienting, whereas the others have at least one case of the intersection property failing. We also note that the only self-dual C-string in Table \ref{table3} is the one with symbol $[3,4,4,3]$.

Of course the more normal subgroups a group has, the more stringent the unravelled condition becomes. We close this section including an example of a soluble group which possess an unravelled C-string. 

\begin{personalEnvironment}\label{(5.5)}$\boldsymbol{\mkern-7mu )}$
$G$ of order $1296 = 2^4.3^4.$
\end{personalEnvironment}
Let $t_1,t_2,t_3,t_4$ be the elements of $\Sym{27}$ as follows:-
\begin{align*}
    t_1=&(4,10)(7,15)(9,17)(12,20)(14,22)(16,23)(19,25)(21,26)(24,27), \\
    t_2 =&(2,4)(5,10)(6,9)(11,17)(12,15)(13,16)(18,23)(19,22)(24,26),\\
    t_3 =&(2,3)(5,8)(7,9)(11,13)(12,16)(15,17)(19,21)(20,23)(25,26) \text{ and }\\
    t_4 =&(1,3)(2,6)(4,9)(5,11)(7,14)(10,17)(12,19)(15,22)(20,25).
\end{align*}
Set $G= \langle t_1,t_2,t_3,t_4\rangle$. Then $\{t_1,t_2,t_3,t_4\}$ is an unravelled C-string for $G$ with diagram \tikzstyle{int}=[draw, fill=black!20, minimum size=2em]
\tikzstyle{init} = [pin edge={to-,thin,black}]
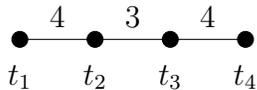
\begin{tikzpicture}[baseline={(0,-0.5)}]

\draw (0,-0.5) -- (1,-0.5);
\draw (1,-0.5) -- (2,-0.5);
\draw (2,-0.5) -- (3,-0.5);
\filldraw [black] (0,-0.5) circle (3pt);
\filldraw [black] (1,-0.5) circle (3pt);
\filldraw [black] (2,-0.5) circle (3pt);
\filldraw [black] (3,-0.5) circle (3pt);
\node (p1) at (0.5,-0.2) {$4$};
\node (p2) at (1.5,-0.2) {$3$};
\node (p3) at (2.5,-0.2) {$4$};

\node (t1) at (0,-1) {$t_1$};
\node (t2) at (1,-1) {$t_2$};
\node (t3) at (2,-1) {$t_3$};
\node (t4) at (3,-1) {$t_4$};
\end{tikzpicture} and Betti numbers $[1,27,81,81,27,1]$. Evidently both vertices and facets of this polytope have the Coxeter group $B_3$ as their automorphism groups. The diameter of the chamber graph is 18, and for $\gamma$ a chamber the disc sizes are as follows 

\begin{table}[H]
    \centering
\begin{tabular}{c | c c c c c c c c c}
    i&1&2&3&4&5&6&7&8&9\\
    \hline
    $|\Delta_i(\gamma)|$   &4 & 9 & 17
& 28 & 42 & 60 & 81 & 105 & 129  \\ 
\end{tabular}\\
\begin{tabular}{c c c c c c c c c }
10&11&12&13&14&15&16&17&18\\
\hline 
 147 & 157 & 155 & 138 & 109 & 71 & 33 & 9 & 1 \\
\end{tabular}
\end{table}

\end{document}